\documentclass[11pt,leqno]{article} 
\usepackage{graphics}

\usepackage{subfigure}
\usepackage{graphicx}

\newtheorem{thm}{Theorem}[section]
\newtheorem{lma}{Lemma}[section]

\newtheorem{prop}{Proposition}

\newcommand{\beqa}{\begin{eqnarray}}
\newcommand{\eeqa}{\end{eqnarray}}

\newcommand{\pf}{\noindent {\bf Proof:} $\s$ }
\newcommand{\epf}{ \hfill$\diamondsuit$ \medskip}

\newcommand{\beq}{\begin{equation}}
\newcommand{\eeq}{\end{equation}}
\newcommand{\lbl}{\label}
\newcommand{\s}{\; \;}

\newcommand{\om}{\overline M}
\newcommand{\um}{\underbar M}

\newcommand{\ra}{\rightarrow}
\newcommand{\al}{\alpha}
\newcommand{\p}{\varphi}

\title{Global solution curves for first order periodic problems, with applications}

\author{
Philip Korman   \\ 
Department of Mathematical Sciences \\ 
University of Cincinnati \\ 
Cincinnati Ohio 45221-0025 \\
\\
Dieter S. Schmidt \\
Department of Electrical Engineering Computer Science \\
University of Cincinnati \\
Cincinnati, Ohio 45221-0030}

\date{}

\begin{document}

\maketitle
\begin{abstract} 
Using  continuation methods and bifurcation theory, we study the exact multiplicity of periodic solutions, and the global solution structure, for periodic problems of first order. The results are applied to a population model with fishing, and to the existence and stability of limit cycles.
We also describe in detail our numerical computations of curves of periodic solutions, and of limit cycles.
 \end{abstract}

\begin{flushleft}
\end{flushleft}

\begin{flushleft}
\end{flushleft}

\section{Introduction}
\setcounter{equation}{0}
\setcounter{thm}{0}
\setcounter{lma}{0}

We study   periodic first order equations
\beq
\label{i1}
u'+g(t,u)=f(t)  \,.
\eeq
Here $g(t+T,u)=g(t,u)$, and $f(t+T)=f(t)$ for all $t \in R$, and we are searching for $T$-periodic solutions, $u(t+T)=u(t)$ for all $t$. Similarly to \cite{K1} and \cite{K2}, we decompose $f(t)=\mu+e(t)$, with $\int _0^T e(t) \, dt=0$ and $\mu \in R$, and $u(t)=\xi+U(t)$, with $\int _0^T U(t) \, dt=0$ and $\xi \in R$, and show the existence of a continuous curve of solutions $(u,\mu)(\xi)$, i.e., $u=u(t,\xi)$, and $\mu=\mu (\xi)$, under rather general assumptions on $g(t,u)$ and $f(t)$. This curve exhausts the set of all periodic solutions of (\ref{i1}). We refer to the continuous function $\mu =\mu (\xi)$ as the {\em solution curve}. To obtain $T$-periodic solutions of 
\beq
\label{i2}
u'+g(t,u)=0  \,,
\eeq
we study the curve of $T$-periodic solutions of    ($\mu \in R$)
\beq
\label{ii3}
u'+g(t,u)=\mu  \,.
\eeq
If the solution curve of (\ref{ii3}), $\mu =\mu (\xi)$, has a root $\xi _0$, then the problem (\ref{i2}) has $T$-periodic solution of average $ \xi _0$. We  give some conditions under which $\mu =\mu (\xi)$ changes sign, and hence it has  roots.
\medskip

As a first application of the above solution curves, we re-derive two results of H. Brezis and L. Nirenberg \cite{BN} on the existence of periodic solutions of (\ref{i2}), providing a wider perspective on these results, and some extra information.   (The paper \cite{BN} dealt with a more general setting of first order equations on a torus.) However, these results have rather stringent conditions, which tend  not to hold in applications. We present further existence results for (\ref{i2}) by studying global properties of the solution curve of (\ref{ii3}), and we also discuss stability properties of the corresponding  $T$-periodic solutions of (\ref{i2}).
\medskip

We apply our results to a periodic population model with fishing
\beq
\lbl{i15}
u'=u(a(t)-u)-\mu f(t) \,,
\eeq
depending on a parameter $\mu>0$. Here $a(t)$ and $f(t)$ are given $T$-periodic functions, and we are looking for a  $T$-periodic  solution $u(t)>0$. We study the curves of  $T$-periodic  solutions of (\ref{i15}) as the  parameter $\mu$ varies. It is known that $T$-period solutions of (\ref{i15}) lie on parabola-like curve with a turn to the left in the $(\mu, ||u||)$ plane, see \cite{A}, \cite{M}, \cite{N}. If $\mu _0$ is the turning point on the solution curve, then $\mu _0 \int_0^t f(t) \, dt$ gives the maximal sustainable level of fishing. If one changes $f(t)$,  $\mu _0$ is changed too, and the natural question is: which $f(t)$ gives the largest  maximal sustainable level of fishing. We prove that in case $a(t)$ is a constant, the optimal $f(t)$ is also a constant. For a general $a(t)$, the optimal $f(t)$ is not necessarily a constant, as our computations show. Based on our numerical experiments,  we state two  open problems.
\medskip

Another application is to limit cycles, that  are closed solutions of autonomous systems of the type
\beqa
\lbl{i10} 
& X'(t)=F(X,Y) \\ \nonumber
& Y'(t)=G(X,Y) \,.\nonumber
\eeqa
It is well known that each limit cycle contains at least one rest point, see e.g., S.H. Strogatz \cite{St}. Let $(x_0,y_0)$ be a rest point of (\ref{i10}), so that $F(x_0,y_0)=G(x_0,y_0)=0$.
We wish to find a  limit cycle, enclosing  $(x_0,y_0)$. Setting $x=X-x_0$ and $y=Y-y_0$, and introducing the polar coordinates $(r,\theta)$ in the $xy$-plane, we obtain a problem of the type
\beq
\label{i3}
r'+g(\theta,r)=0 \,,
\eeq
with $g(\theta+2 \pi,r)=g(\theta,r)$, and we are looking for a $2\pi$-periodic solution $r=r(\theta)$. The existence results of  \cite{BN} relied on a priori estimates. These  estimates were obtained under the following two conditions: either $g_r(\theta,r)>0$ for all $\theta$ and $r$, or the problem (\ref{i10}) has an ordered sub and super solution pair. Both conditions are way too stringent for limit cycles.
The term ``a priori estimate" refers to estimating solutions, without knowledge of these solutions. If one attempts to get a priori estimates also without knowledge of the equation, one is led to impose very stringent conditions, as above. In our existence results we assume that a priori estimates  hold, and then we derive such estimates for concrete systems by using trapping regions, or energy considerations.
Our approach to limit cycles provides a straightforward alternative to the traditional one, based on the Poincare-Bendixson theorem, making it unnecessary to search for the trapping region in each particular case. Moreover, it unifies some diverse applications, and makes it convenient to handle multiple limit cycles. 
\medskip

We also present a detailed description of the numerical computation of the global solution curves, and of limit cycles. In particular, we compute  limit cycles for the  Van der Pol and Sel'kov's equations, and the bifurcation diagrams for periodic population models with ``fishing".

\section{The global solution curve}
\setcounter{equation}{0}
\setcounter{thm}{0}
\setcounter{lma}{0}

We begin with some preliminary results. The following  lemma is well-known as Wirtinger's inequality. Its proof follows easily by using the complex Fourier series, and the orthogonality of the functions $\{e^{i \omega n t} \}$ on the interval $(0,T)$, see e.g., P. Korman \cite{K}.
\begin{lma}\lbl{lma:g1}
Assume that $f(t)$ is a continuously differentiable function  of period $T$, and of zero average, i.e.,  $\int_0^T f(s) \, ds=0$. Then
\[
\int _0^T {f'}^2(t) \, dt \geq \omega ^2 \int _0^T f^2(t) \, dt \,, \s \mbox{with $\omega=\frac{2 \pi}{T}$} \,.
\]
\end{lma}

We shall need the following standard lemma, see e.g., \cite{N}, \cite{K}.

\begin{lma}\lbl{lma:g2}
Consider the equation
\[
z'+a (t)z =0 \,,
\]
with a given continuous $T$-periodic function $a (t)$. This equation has a non-zero  $T$-periodic solution if and only if $\int_0^T a(t) \, dt=0$. 
\end{lma}

The following lemma will be  crucial for continuation of solutions.
\begin{lma}\lbl{lma:g33}
Let $a(t)$ be a given continuous  $T$-periodic function, and $b\in R$ a constant.  The differential equation
\begin{equation}  w'+a(t)w=b  \label{deq1}
\end{equation}
does not have a $T$-periodic solution, $w(t)=w(t+T)$,  satisfying
\beq
\lbl{int}
\int_0^T w(t) \, dt=0 \,,
\eeq
unless $b=0$, in which case $w(t)=0$.
\end{lma}

\pf
With $\nu(t)=e^{\int_0^t a(\tau)  \, d\tau} >0$, the solution to (\ref{deq1}) satisfying $w(0)=w_0$ can be written as
\[
w(t)=\frac{1}{\nu(t)} \left( b\int_0^t \nu(s)  \, ds  +w_0 \right) \,.
\]
The periodicity condition  $w(T)-w(0)=0$, and the requirement (\ref{int})  lead to the following two equations 
\begin{eqnarray*}
b\int_0^T \nu(s) ds + w_0(1-\nu(T))&=&0,\\
b\int_0^T\frac{\int_0^t\nu(s)  \, ds}{\nu(t)}\; dt +w_0\int_0^T \frac{1}{\nu(t)}  \, dt&=&0.
\end{eqnarray*}
The determinant of this linear system for $b$ and $w_0$ is 
\[
D=\int_0^T \nu(s)  \,  ds \int_0^T\frac{dt}{\nu(t)} \, dt  -\int^T_0\frac{\int_0^t \nu(s)  \,  ds}{\nu(t)}  \,  dt+\nu(T)\int_0^T\frac{\int_0^t \nu(s)  \, ds}{\nu(t)}  \, dt.
\]
Since $\nu(t)>0$ it follows that $\int_0^t\nu(s) \, ds < \int_0^T \nu(s) \, ds$ for $t<T$, and then the sum of the first two terms  in $D$ is positive.  This shows that $D>0$, and then the linear system has  only the trivial solution: $b=0$ and $w_0=0$.
\epf

We shall denote by $C_T$ and $C^1_T$ the subspaces of $C(R)$ and $C^1(R)$ respectively, consisting of $T$-periodic functions. The next lemma represents the well-known {\em Fredholm alternative}. It can be proved by direct integration.
\begin{lma}\lbl{lma:g4}
Consider the equation
\beq
\label{g77}
y'+a(t)y=f(t) \,,
\eeq
with $a(t)$, $f(t) \in C_T$. If $\int _0^T a(t) \, dt \ne 0$, the equation (\ref{g77}) has a unique solution $y(t) \in C^1_T$, for any  $f(t) \in C_T$.
In case $\int _0^T a(t) \, dt = 0$, the equation (\ref{g77}) has $T$-periodic solutions if and only if $\int _0^T f(t) u_0(t) \, dt = 0$, where $u_0(t) \in C_T$ is any non-trivial solution of 
\beq
\label{g7a}
u'-a(t)u=0 \,.
\eeq
The number of  $T$-periodic solutions of (\ref{g77}) is then infinite.
\end{lma}

The following lemma will be used to continue solutions of fixed average. 

\begin{lma}\lbl{lma:g5}
Let $a(t) \in C_T$.
For any $\xi \in R$ and $e(t)  \in C_T$, there exists a pair  $(\mu ^*,z(t))$, with $\mu ^* \in R$ and $z(t)  \in C^1_T$, solving
\beqa
\label{g88}
& z'+a(t)z=\mu ^*+e(t) \\ \nonumber
& \frac{1}{T} \int _0^T z(t) \, dt =\xi \,. \nonumber
\eeqa
\end{lma}

\pf
Case 1. $\int _0^T a(t) \, dt \ne 0$. Consider the operator $L \,: C_T^1 \ra C_T$ given by $L(z)=z'+a(t)z$. By Lemma \ref{lma:g4}, $L^{-1}$ exists, and we can express the solution of the first equation in (\ref{g88}) as 
\[
z=\mu ^* L^{-1}(1)+L^{-1}(e(t)) \,.
\]
By Lemma \ref{lma:g33}, $\int _0^T L^{-1}(1) \, dt \ne 0$, and hence we can choose $\mu ^*$ to satisfy the second equation in (\ref{g88}).
\medskip

\noindent
Case 2. $\int _0^T a(t) \, dt = 0$. This time $L$ is not invertible. Let $u_0(t)>0$ be any non-trivial solution of (\ref{g7a}). Choose $\mu ^*$ so that $\int _0^T \left[\mu ^*+e(t) \right] u_0(t)\, dt = 0$. Then  the first equation in (\ref{g88}) has infinitely many solutions of the form $z(t)=z_0(t)+c \bar u_0(t)$, where $z_0(t)$ is some fixed solution, and $\bar u_0(t)>0$ is any $T$-periodic solution of
\[
z'+a(t)z=0 \,.
\]
Choose the constant $c$ so that  the second equation in (\ref{g88}) is satisfied.
\epf

\begin{lma}\lbl{lma:sub}
Let $u(t) \in C^1_T$ be   a solution of 
\beq
\lbl{sub1}
u'(t)=f(t,u) \s\s \mbox{for all $t \in R$} \,,
\eeq
for some continuous $f(t,u)$, while $v(t) \in C^1_T$ is a sub (super) solution of (\ref{sub1}), so that
\[
v'(t)<f(t,v) \s \left( v'(t)>f(t,v) \right) \s\s \mbox{for all $t \in R$} \,.
\]
Then $u(t)$ and $v(t)$ do not intersect (and hence they are strictly ordered).
\end{lma}

\pf
Assume that $v(t)$ is a super  solution, and the other case is similar. If $t_0$ is a point of intersection of  $u(t)$ and $v(t)$, then $v'(t_0)>u'(t_0)$.
If the $T$-periodic functions  $u(t)$ and $v(t)$ intersected, they would intersect at least twice, and at one of the intersection points we would have $v'(t_0) \leq u'(t_0)$, a contradiction.
\epf

Any $f(t) \in C_T$ can be written in the form $f(t)=\mu +e(t)$, where $\mu \in R$, and $e(t)\in C_T$  satisfies $\int _0^T e(t) \, dt=0$.
We now consider the following problem: find $u(t) \in C^1_T$ and $\mu \in R$ satisfying
\beq
\label{g3}
u'+g(t,u)=\mu +e(t) \,, \s\s \mbox{$u(t+T)=u(t)$} \; \mbox{ for all $t$} \,.
\eeq

\begin{thm}\lbl{thm:5}
Assume that $g(t,u) \in C^{0,1}(R \times R)$ satisfies $g(t+T,u)=g(t,u)$ for all $t$ and $u$, and $e(t)\in C_T$. Assume that either $g(t,u)=g(u)$, or there is a constant $c_0>0$ such that 
\beq
\label{g4}
|g(t,u)| \leq c_0 \,, \s\s \mbox{ for all $t \in R$, and $u\in R$} \,.
\eeq
Then for any $\xi \in R$, one can find a unique solution pair $(u(t),\mu)$ of (\ref{g3}), with $\frac{1}{T} \int _0^T u(t) \, dt=\xi$. Moreover, the solution $u(t,\xi)$ is increasing in $\xi$.
\end{thm}

\pf
Decompose $u(t)=\xi+U(t)$, with $\xi \in R$, and $\int _0^T U(t) \, dt=0$. Then $U(t)$ satisfies
\beq
\label{g5}
U'+g(t,\xi+U)=\mu +e(t) \,, \s \mbox{$U(t+T)=U(t)$} \; \mbox{  for all $t$} \,.
\eeq
We show next that $\max _t |U(t)|$ is bounded, uniformly in $\xi$ and $\mu$.
Multiply (\ref{g5}) by $U'(t)$ and integrate over $(0,T)$. Obtain
\beq
\label{g6}
\int _0^T {U'}^2 \, dt +\int _0^T g(t,\xi+U) U' \, dt=\int _0^T U' e(t) \, dt \,.
\eeq
If $g(t,u)=g(u)$, the second integral on the left vanishes by the periodicity of $U(t)$, and 
 we get an estimate of $\int_0^T {U'}^2(t) \,dt$. In case $g(t,\xi+U)$ is bounded, we estimate the second integral on the left by
$\epsilon   \int_0^T {U'}^2(t) \,dt+c_1(\epsilon)$, choose $\epsilon$ small, and again obtain  an estimate of $\int_0^T {U'}^2(t) \,dt$. 
 We then estimate $\int_0^T {U}^2(t) \,dt$ by Wirtinger's inequality, and conclude a bound on  $\max |U(t)|$ by Sobolev's embedding. 
\medskip

We now embed (\ref{g3}) into a family of problems
\beq
\label{g7}
u'+kg(t,u)=\mu +e(t) \,, \s \mbox{$u(t+T)=u(t)$}  \; \mbox{ for all $t$} \,.
\eeq
depending on a parameter $k$, with $0 \leq k \leq 1$ ($k=1$ corresponds to  (\ref{g3})).  When $k=0$, and $\mu=0$, the problem has infinitely many solutions in $C^1_T$, so that for any $\xi ^0 \in R$ we can choose the unique solution $u_0(t)$ with $\frac{1}{T} \int _0^T u_0(t) \, dt=\xi ^0$.  We now continue in $k$ the solutions of
\beqa
\label{g8}
& F(u,\mu,k) \equiv u'+kg(t,u)-\mu -e(t)=0 \,, \s \mbox{$u(t+T)=u(t)$} \\ \nonumber
& \frac{1}{T} \int _0^T u(t) \, dt=\xi ^0 \,, \nonumber
\eeqa
with the operator $F(u,\mu,k) \, : C^1_T \times R \times R \ra C_T$. We  show next that the implicit function theorem applies, allowing us to continue the solution pair $(u,\mu)$ as a function of $k$. Compute the Frechet derivative (with $w \in C^1_T$)
\beqa \nonumber
& F_{(u,\mu)}(u,\mu, k)(w, \mu^*)= w'+kg_u(t,u)w-\mu^*  \,,  \\ \nonumber
& \int _0^T w(t) \, dt=0 \,. \nonumber
\eeqa
By Lemma \ref{lma:g33}, the map $F_{(u,\mu)}(u,\mu, k)(w, \mu^*)$ is injective, and by Lemma \ref{lma:g5} this map is surjective. Hence, the implicit function theorem applies, and we have a solution curve $(u,\mu)(k)$. By the a priori estimate above, this curve continues for all $0 \leq k \leq 1$, and at $k=1$ we obtain a solution of the problem (\ref{g3}), with $\frac{1}{T} \int _0^T u(t) \, dt=\xi ^0$.
\medskip

Turning to the uniqueness part, let $(\bar u(t),\bar \mu)$ be another solution of (\ref{g3}),  with $\frac{1}{T} \int _0^T \bar u(t) \, dt=\xi ^0$. Then $(\bar u(t),\bar \mu)$ is a solution of (\ref{g7}) at $k=1$. We continue this solution backward in $k$, until $k=0$, using  the implicit function  theorem. Integrating the equation, we conclude that  $\mu =0$, when $k=0$. Then at  $k=0$, $\mu =0$ we obtain a  $T$-periodic solution of $u'=e(t)$, which is different from $u_0(t)$ (since the solution curves do not intersect, by the uniqueness part of the implicit function  theorem), a contradiction.
\medskip

Next, we use  the implicit function theorem to show that  solutions of (\ref{g3}) can be continued in $\xi$, for all $\xi \in R$. Decomposing $u(t)=\xi+U(t)$ with $\int _0^T U(t) \, dt=0$, we see that $U(t)$ satisfies
\beqa \nonumber
& F(U,\mu,\xi) \equiv U'+kg(t,\xi+U)-\mu -e(t)=0 \,, \s \mbox{$U(t+T)=U(t)$} \\ \nonumber
& \int _0^T U(t) \, dt=0 \,, \nonumber
\eeqa
with the operator $F(U,\mu,\xi) \, : C^1_T \times R \times R \ra C_T$. Compute the Frechet derivative (with $w \in C^1_T$)
\beqa \nonumber
& F_{(U,\mu)}(U,\mu, k)(w, \mu^*)= w'+kg_u(t,u)w-\mu^*  \,,  \\ \nonumber
& \int _0^T w(t) \, dt=0 \,. \nonumber
\eeqa
As above, we see that  the implicit function theorem applies, and we have a smooth solution curve $(u,\mu)(\xi)$ for the problem  (\ref{g3}). By the a priori estimate above, this curve continues for all $\xi \in R$.
\medskip

Finally, we prove the monotonicity of solutions. Let $\xi _1<\xi _2$ be arbitrary. Assume, first, that $\mu(\xi _2)=\mu (\xi _1)$. Then $u(t,\xi _1)$ and $u(t,\xi _2)$ are two solutions of the same equation, and they cannot intersect by uniqueness of the initial value problems. Since  $u(t,\xi _2)$ has a greater average, it follows that $u(t,\xi _2)>u(t,\xi _1)$ for all $t$. Assume next that $\mu(\xi _2)<\mu (\xi _1)$ ($\mu(\xi _2)>\mu (\xi _1)$). Then $u(t,\xi _2)$ is sub (super) solution for the equation that  $u(t,\xi _1)$ satisfies. By Lemma \ref{lma:sub}, $u(t,\xi _1)$ and $u(t,\xi _2)$ are ordered, and hence $u(t,\xi _2)>u(t,\xi _1)$ for all $t$.
\epf

We now discuss the stability of $T$-periodic solutions of
\beq
\lbl{sta1}
u'+g(t,u)=e(t) \,,
\eeq
with $g(t,u)$ and $e(t)$ satisfying $\int_0^T e(t) \, dt=0$, and the conditions of the Theorem \ref{thm:5} above.

\begin{thm}\lbl{thm:stab}
Let $\mu=\mu (\xi)$ be the solution curve of (\ref{g3}). Assume that $\mu (\xi _0)=0$, so that (\ref{sta1}) has a solution of the form $u_0(t)=\xi _0+U(t)$ with $\int_0^T U(t) \, dt=0$. Assume that there is an interval $(\xi _0-\epsilon,\xi _0+\epsilon)$ so that $\mu (\xi)<0$ ($\mu (\xi)>0$) on $(\xi _0-\epsilon,\xi _0)$, and $\mu (\xi)>0$ ($\mu (\xi)<0$)  on $(\xi _0,\xi _0+\epsilon)$, for some $\epsilon >0$. Then $u(t)$ is a stable (unstable)  solution of  (\ref{sta1}).
\end{thm}

\pf
Assume that the first set of inequalities holds, and the other case is similar.
Take any $\xi _1 \in (\xi _0-\epsilon,\xi _0)$, calculate $\mu _1=\mu (\xi _1)<0$, and the corresponding solution $u_1(t)$  of (\ref{g3}). Then
\beq
\lbl{20}
u_1'+g(t,u_1)=\mu _1+e(t)<e(t) \,,
\eeq
so that $u_1(t)$ is a sub solution of  (\ref{sta1}). By Lemma \ref{lma:sub}
\beq
\lbl{21}
u_1(t)< u_0(t) \s \mbox{for all $t$} \,.
\eeq

Thus we have a family of sub solutions below $u_0 (t)$, and tending to $u_0 (t)$. Similarly, we produce a family of super solutions above $u_0 (t)$, and tending to $u_0 (t)$. The stability of $u_0 (t)$ is then well-known. (By Lemma \ref{lma:sub}, any solution of (\ref{sta1}), lying below $u_0(t)$, at any point cuts above the subsolution $u_1(t)$, passing through the same point.)
\epf

\section{Existence of $T$-periodic solutions}
\setcounter{equation}{0}
\setcounter{thm}{0}
\setcounter{lma}{0}
We apply the global solution curves to give alternative proofs of two results of  H. Brezis and L. Nirenberg \cite{BN} on the existence of $T$-periodic solutions of 
\beq
\lbl{1}
u'+g(t,u)=0 \,, \; u(t)=u(t+T)  \;\; \mbox{ for all $t$} \,.
\eeq

\begin{thm}\lbl{thm:1}
Assume that $g(t,u) \in C^{0,1}(R \times R)$ satisfies $g(t+T,u)=g(t,u)$ for all $t$ and $u$. Assume the existence of  constants $\underbar M<\overline M$ such that
\beq
\lbl{2}
\int_0^T g(t,\underbar M) \, dt<0 \,, \s \int_0^T g(t,\om) \, dt>0 \,.
\eeq
Assume finally the  existence of  constants $M_1<M_2$  such that for all $t$
\beqa
\lbl{3}
& g(t,u)>g(t,\om) \s \mbox{for $u>M_2$} \\
& g(t,u)<g(t,\um) \s \mbox{for $u<M_1$} \,. \nonumber
\eeqa
Then the equation (\ref{1}) has a $T$-periodic solution.
\end{thm}

\pf
Step 1. A priori estimate. Write $g(t,\om)=\xi -v'(t)$, with $\int_0^T v(t) \, dt=0$, and $\xi>0$ by (\ref{2}). Consider $u(t)-v(t)$. It satisfies
\beq
\lbl{3b}
\left(u-v \right)'+g(t,u)=-v'=g(t,\om)-\xi<g(t,\om) \,.
\eeq
Let $t_0$ be the point of maximum of $u-v$. At $t_0$
\beq
\lbl{3c}
g(t_0,u(t_0))<g(t_0,\om) \,,
\eeq
so that $u(t_0)<M_2$ by (\ref{3}). Then
\[
u-v \leq \max \left(u-v \right)=u(t_0)-v(t_0) <M_2+ \max |v| \,,
\]
\[
u <v+M_2+ \max |v| \leq M_2+2 \max |v| \,.
\]
This gives  an estimate of $u(t)$ from above. To get an estimate from below, write $g(t,\underline M)=\eta+w'$, with $\int_0^T w(t) \, dt=0$, and $\eta < 0$ by (\ref{2}). We have
\[
\left(u+w \right)'+g(t,u)=w'=g(t,\underline M)-\eta>g(t,\underline M ) \,.
\]
Let $t_1$ be the  point of minimum of $u+w$. Then $g(t_1,u(t_1))>g(t_1,\underline M )$, and $u(t_1)>M_1$ by (\ref{3}). We conclude
\[
u+w \geq \min \left(u+w \right)=u(t_1)+w(t_1)>M_1-\max |w| \,,
\]
\[
u>-w+M_1-\max |w|\geq M_1-2\max |w| \,.
\]

\noindent
Step 2. Modify $g(t,u)$ for $u$ outside of $[\underline M, \om]$ to be a bounded function, satisfying (\ref{3}), and call it $\bar g(t,u)$. Clearly, $\bar g(t,u)$ satisfies (\ref{2}) too. By above, any $T$-periodic solution of 
\[
u'+\bar g(t,u)=0
\]
satisfies the same a priori estimate as before, and hence it is a solution of (\ref{1}). Therefore we may still work with  (\ref{1}), but assume that $g(t,u)$ is bounded.
\medskip

\noindent
Step 3. Embed  (\ref{1}) into a family of problems
\beq
\lbl{5}
u'+g(t,u)=\mu \,, \; u(t)=u(t+T) \;\; \mbox{ for all $t$} \,,
\eeq
depending on a parameter $\mu$. Decompose $u(t)=\xi+U(t)$, with $\xi \in R$ and $\int_0^T U(t) \,dt=0$. By Theorem \ref{thm:5} we have a continuous solution curve $\mu=\mu(\xi)$ for all $\xi \in (-\infty,\infty)$. Recall that $U(t)$ satisfies 
\beq
\lbl{6}
U'+g(t,\xi+U)=\mu \,, \; \int_0^T U(t) \,dt=0 \,, \; U(t)=U(t+T) \,,
\eeq
and $\max |U(t)|$ is bounded uniformly in $\xi$ and $\mu$.

\noindent
Step 4. Integrate (\ref{6}):
\[
\mu T=   \int_0^T g(t,\xi+U(t)) \,dt  \,.
\]
For $\xi$ large, $\xi+U(t)>M_2$ for all $t$. Using  (\ref{3}), and then  (\ref{2})
\beq
\lbl{6a}
\mu T>   \int_0^T g(t,\om) \,dt  >0\,,
\eeq
so that $\mu >0$. Similarly, $\mu <0$ for $-\xi$ large. We conclude the existence of $\xi _0$ with $\mu (\xi _0)=0$, which gives a solution of (\ref{1}).
\epf

\noindent
{\bf Remarks} 
\begin{enumerate}
  \item  Schauder's fixed point theorem was used in \cite{BN}, and  as mentioned in \cite{BN}, this result can also be proved using super and sub solutions.
  \item If $g(t,0)=0$ for all $t$, and $g(t,u)$ is Lipschitz continuous in $u$, then either $u(t)>0$ or  $u(t)<0$  for all $t$, by uniqueness for initial value problems.
  \item  The first condition (\ref{3}) is rather restrictive. In particular, it does not apply to the Sel'kov system, considered below.
  \item Our approach provides some extra information: there exists a curve $(\mu, u(t))(\xi)$ of $T$-periodic solutions of (\ref{5}), and $\xi$ (the average value of $u(t)$) is the global parameter on this curve.
\end{enumerate}  
\medskip

The following variation on this result is also due to H. Brezis and L. Nirenberg \cite{BN}.
\begin{thm}\lbl{thm:3}
Assume that $g(t,u) \in C^{0,1}(R \times R)$ satisfies $g(t+T,u)=g(t,u)$, and the condition (\ref{2}).
Assume also that 
\[
g_u(t,u)>0\s \mbox{for all $t$,$\, u$} \,. \nonumber
\]
Then the equation (\ref{1}) has a $T$-periodic solution.
\end{thm}

\pf
The argument is almost the same. From (\ref{3c}), it follows that $u(t_0)<\bar M$, leading to the a priori estimate as above. As for the estimate (\ref{6a}), once $\xi$ is so large that $\xi+U(t)> \overline M$, we have 
\[
\mu T=\int_0^T g(t,\xi+U(t)) \,dt>   \int_0^T g(t,\overline M) \,dt  >0\,,
\]
and the rest of the proof is the same.
\epf

\section{A further existence result, and  applications to limit cycles}
\setcounter{equation}{0}
\setcounter{thm}{0}
\setcounter{lma}{0}

We now change the notation, and consider  the existence of $2\pi$-periodic solutions of  
\beq
\lbl{12}
\frac{dr}{d \theta}+g(\theta,r)=0 \,,
\eeq
where $g(\theta,r)$
is $2\pi$-periodic in $\theta$. It is known that any  $2 \pi$-periodic solution of (\ref{12}) is in fact a limit cycle, attracting the nearby solutions, see V.A. Pliss \cite{P}. We shall assume that
\[
g(\theta,0)=0 \,, \s \mbox{for all $\theta \in [0,2\pi)$} \,,
\]
which is the case in  applications to limit cycles. The last assumption implies that solutions of (\ref{12}) do not change sign, and we shall search for solutions with $r(\theta)>0$ for all $\theta$. We now give a general result on the existence of a $2\pi$-periodic solution of (\ref{12}).

\begin{thm}\lbl{thm:2}
Assume that $g(\theta,0)\equiv 0$, and $g(\theta,r) \in C^{0,1}(\bar R_+ \times [0,2 \pi])$ satisfies $g(\theta+2 \pi,r)=g(\theta,r)$, for all $\theta$ and $r$. Assume that there a constant $M_0>0$, such that any positive $2 \pi$-periodic solution of (\ref{12}) satisfies $0<r(\theta)<M_0$. Assume that for some $p \geq 1$ the following limits exist, uniformly in $\theta$
\beq
\lbl{12a}
\lim _{r \ra 0} \frac{g(\theta,r)}{r}=a(\theta) \,, \s \mbox{with $\int_0^{2 \pi} a(\theta) \, d\theta >0$ $(<0)$} \,,
\eeq
\beq
\lbl{12b}
\lim _{r \ra \infty} \frac{g(\theta,r)}{r^p}=b(\theta) \,, \s \mbox{with $\int_0^{2 \pi} b(\theta) \, d\theta <0$ $(>0)$} \,.
\eeq
Then the equation  (\ref{12}) has  a $2\pi$-periodic solution $r(\theta)>0$.
\end{thm}

\pf
Assume that the first set of inequalities holds in (\ref{12a}), (\ref{12b}), and the other case is similar.
By modifying $g(\theta,r)$ for $r>M_0$, and for $r<0$, we may assume that this function is bounded, for all $\theta$ and $r$.
We embed (\ref{12}) into a family of problems
\beq
\lbl{12c}
r'+g(\theta,r)=\mu \,, \s r(\theta+2\pi )=r(\theta)\,.
\eeq
Decompose $r(\theta)=\xi+R(\theta)$, with $\int_0^{2\pi} R(\theta) \, d \theta=0$. By Theorem \ref{thm:5}, we have a solution curve $\mu=\mu (\xi)$, with $\xi \geq 0$, and $\mu (0)=0$. We wish to show that $\mu (\xi _0)=0$ for some $\xi _0 >0$.
\medskip

We claim that for any $\xi >0$, $r(\theta)>0$ for all $\theta$. In case $\mu=0$, $r(\theta)$ cannot change sign by the uniqueness theorem, and since its average is positive, $r(\theta)$ is positive. In case $|\mu|>0$, we see from the equation (\ref{12c}) that $r'(\theta)$ has the same sign at all of its roots, which makes roots impossible. Again, $r(\theta)$ is positive, since its average $ \xi$ is positive.
\medskip

We show next that $\mu (\xi)>0$ for $\xi>0$ and small.  Divide  the equation (\ref{12c})  by $r(\theta)>0$, and integrate over $(0,2 \pi)$. By periodicity, 
$\int_0^{2\pi} \frac{r'}{r} d \theta=0$, and we have 
\[
\int_0^{2\pi} \frac{g(\theta,r(\theta))}{r(\theta)} \, d \theta=\mu \int_0^{2\pi} \frac{1}{r(\theta)} \, d \theta\,.
\]
By continuity, $r(\theta)$ is small when $\xi>0$ and small. Then $\mu >0$, in view of the condition (\ref{12a}).
Similarly, dividing  the equation (\ref{12c})  by $r^p(\theta)>0$, and integrating  over $(0,2 \pi)$, we see that $\mu(\xi)>0$ for $\xi$ large. By continuity of $\mu(\xi)$, we conclude that $\mu (\xi _0)=0$ at some $\xi _0 >0$.
\epf

\noindent
{\bf Example} This theorem provides the existence of positive $2\pi$-periodic solution  to the equation
\[
r'+a(\theta)r+b(\theta)r^p=0 \,,
\]
with $p>1$, and $2\pi$-periodic $a(\theta)$ and $b(\theta)$ such that $\int_0^{2 \pi} a(\theta) \, d\theta $ and $\int_0^{2 \pi} b(\theta) \, d\theta $ have opposite signs, provided there is an a priori bound on $2\pi$-periodic solutions. An   a priori bound will hold if, for example, $b(\theta)$ has no roots. Both Theorems \ref{thm:1} and \ref{thm:3} do not apply in case of sign-changing $a(\theta)$ and $b(\theta)$.
\medskip

We now provide a version of the last result. Introduce the function ($M$ is a real number)
\[
\p (M)=\int_0^{2 \pi} g(\theta,M) \, d \theta \,.
\]
Then the conditions (\ref{12a}) and (\ref{12b}) can be replaced by $\p'(0)>0$ ($ <0$), and $\lim _{M \ra \infty} \frac{\p (M)}{M^p}<0$ ($>0$) respectively. Indeed, 
\[
\int_0^{2 \pi} a(\theta) \, d\theta=\int_0^{2 \pi} \lim _{r \ra 0} \frac{g(\theta,r)}{r} \, d \theta=\int_0^{2 \pi} g_r(\theta,0) \, d \theta=\p '(0) \,,
\]
\[
\int_0^{2 \pi} b(\theta) \, d\theta=\int_0^{2 \pi} \lim _{r \ra \infty} \frac{g(\theta,r)}{r^p} \, d \theta=\lim _{M \ra \infty} \frac{\p (M)}{M^p} \,.
\]
\medskip

We now turn to the limit cycles, which are closed solutions of autonomous systems of the type
\beqa
\lbl{10} 
& X'(t)=F(X,Y) \\ \nonumber
& Y'(t)=G(X,Y) \,.\nonumber
\eeqa
It is well known that each limit cycle contains at least one rest point, see e.g., S.H. Strogatz \cite{St}. Let $(x_0,y_0)$ be a rest point of (\ref{10}), so that
\[
F(x_0,y_0)=G(x_0,y_0)=0 \,.
\]
We wish to find a  limit cycle, enclosing  $(x_0,y_0)$. Setting $x=X-x_0$ and $y=Y-y_0$, we obtain
\beqa
\lbl{11} 
& x'(t)=F(x+x_0,y+y_0) \\ \nonumber
& y'(t)=G(x+x_0,y+y_0) \,,\nonumber
\eeqa
with the rest point at the origin $(0,0)$. We  then switch to  polar coordinates in the $(x,y)$ plane, given by $x=r \cos \theta$ and $y=r \sin \theta$. We obtain
\beq
\lbl{11a}
\frac{dr}{dt}=\frac{xF(x+x_0,y+y_0)+y G(x+x_0,y+y_0)}{\sqrt{x^2+y^2}} \,,
\eeq
and
\[
\frac{d \theta}{dt}=\frac{\frac{dy}{dt} x-\frac{dx}{dt} y}{x^2+y^2}=\frac{xG(x+x_0,y+y_0)-y F(x+x_0,y+y_0)}{x^2+y^2} \,.
\]
Using that $\frac{dr}{dt}=\frac{dr}{d \theta}\frac{d \theta}{dt}$, and 
 defining
\[
h(x,y)=\frac{xF(x+x_0,y+y_0)+y G(x+x_0,y+y_0)}{xG(x+x_0,y+y_0)-y F(x+x_0,y+y_0)} \sqrt{x^2+y^2} \,,
\]
and 
\beq
\lbl{17-a}
g(\theta,r)=- h\left (r \cos \theta,r \sin \theta \right) \,,
\eeq
we obtain from (\ref{11a})
\[
\frac{dr}{d \theta}+g(\theta,r)=0 \,.
\]
This is a particular case of  the equation (\ref{12}),
where $g(\theta,r)$ is given by (\ref{17-a}).
Clearly, $g(\theta,0)=0$, and $g(\theta,r)$ has period $2\pi$ in $\theta$. We  study the existence of a $2\pi$-periodic solution $r(\theta)$ of that equation, which corresponds to a limit cycle enclosing  $(x_0,y_0)$, that  is a graph of a polar curve around  $(x_0,y_0)$. Observe that the function  $g(\theta,r)$ is {\em singular}, i.e., it assumes infinite values at some finite $\theta$ and $r$. Singular problems occur often in applications, see the recent book of P.J. Torres \cite{T} for a nice review, or P. Korman \cite{K9}.
\medskip

\noindent
{\bf Example} We consider next the classical Van der Pol equation
\beq
\lbl{17}
x''+(x^2-1)x'+x=0 \,.
\eeq
Setting $x'=y$, we transform it to the system (\ref{10}), where $F(x,y)=y$, $G(x,y)=-(x^2-1)y-x$, with the rest point $x_0=y_0=0$.  To get an a priori bound from above for any limit cycle of the equation (\ref{17}), one can consider the energy, which is decreasing for $|x|>1$.  We search for  a limit cycle, which  is a polar curve in the $(x,y)$ plane, so that  $r=r(\theta)$ satisfies (\ref{12}), with the function in (\ref{17-a}) given by 
\beq
\lbl{l8}
 g(\theta,r)=\frac{2 r (1- \cos{2\theta})(2-r^2(1+\cos{2\theta}))}{8-2(2-r^2)\sin{2\theta}+r^2\sin{4\theta}}.
\eeq

\begin{figure}[ht!]

\subfigure[]{\includegraphics[scale=0.6]{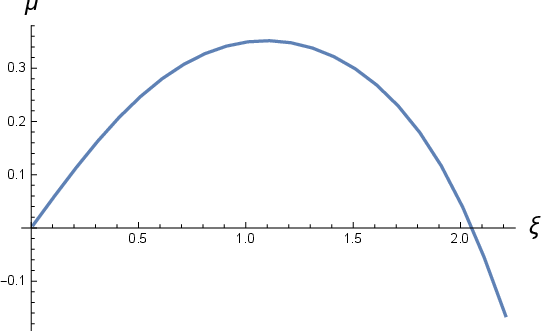}}
\hfill
\subfigure[]{\includegraphics[scale=0.7]{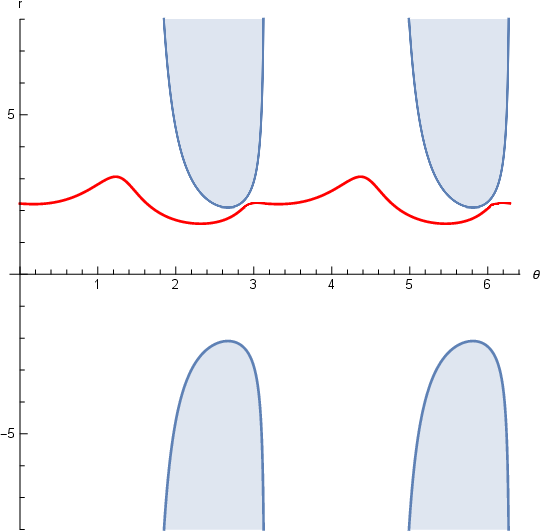}}
\caption{ a)\ The curve of $2\pi$-periodic solutions of (\ref{12c}) for the Van der Pol equation.  b)\ The domain for the $2\pi$-periodic solutions, with a periodic solution with $\xi=2.2$ }
\label{fig:1}
\end{figure}

 To find a $2\pi$-periodic solution of  (\ref{12}), we compute the curve $\mu=\mu (\xi)$  of  $2\pi$-periodic solutions of   (\ref{12c}), and the result is presented in Figure \ref{fig:1}a. One sees that for a $\xi _0$ slightly over $2$, $\mu (\xi _0)=0$, corresponding to a limit cycle  of the  Van der Pol equation. (The value of $\xi _0$ gives the average of $r(\theta)$.) We confirmed by a direct integration that the solution $r=r(\theta)$ at $\mu =0$ is identical to the limit cycle of the  Van der Pol equation. 
\medskip

However, the function $g(\theta ,r)$ is unbounded when the denominator in (\ref{l8}) is zero.  Therefore, the $2\pi$-periodic solutions of (\ref{12c}) will be restricted to the domain defined by $8-2(2-r^2)\sin{2\theta}+r^2\sin{4\theta}=0$, which is displayed in Figure \ref{fig:1}b together with a solution $r(\theta)$ for $\xi=2.2$. With increasing $\xi$ the solutions come even closer to the boundary of this domain and the periodic solutions cease to exist for $\xi>\xi_1$ with $\xi_1 \in (2.2,2.3)$.  It also means that the curve $\mu(\xi)$ can't be continued past $\xi_1$.

\medskip

\noindent
{\bf Example} To model oscillations connected to {\em glycolysis} (the process of cells breaking down sugar to obtain energy) the following model was proposed by E.E. Sel'kov \cite{Sel}, see S.H. Strogatz \cite{St} for a nice presentation,
\beqa
\lbl{be3}
& x'=-x+ay+x^2y \\ \nonumber
& y'=b-ay-x^2y \,.
\eeqa
Here $a$ and $b$ are positive parameters. The unknown functions $x(t)$ and $y(t)$ represent some biological quantities, and are also assumed to be positive.
\medskip

One checks that the only rest point is at $x_0=b$, $y_0=\frac{b}{a+b^2}$. 
It is known that there is a limit cycle around $(x_0,y_0)$, provided that  the parameters $a$ and $b$ satisfy Sel'kov's condition
\beq
\lbl{15}
a^2+2 a
   b^2+ a+ b^4-
   b^2 <0 \,.
\eeq
(This condition holds for small $b>0$, and even smaller $a>0$.)  Next we derive an a priori bound from above for any $2\pi$ periodic solution of the equation (\ref{12}), which corresponds to (\ref{be3}). Since the limit cycle is around $(x_0>0,y_0>0)$, at least part of it must lie in the first quadrant. However, from the equations in (\ref{be3}), the limit cycle  cannot exit  the first quadrant, and hence it lies completely in  the first quadrant.  Adding the equations in (\ref{be3}), we get
\[
\left(x+y \right)'=b-x \,,
\]
so that  in the first quadrant all  trajectories travel southwest if $x>b$. Also, $y'<0$ for $y>\frac{b}{a}$. It follows that  all  trajectories  in the first quadrant enter the {\em trapping region} bounded by the line $y=\frac{b}{a}$ for $0<x \leq b$, the line $x+y=b+\frac{b}{a}$ for $b<x<b+\frac{b}{a}$, and the coordinate axes. In particular, any  limit cycle must lie in this region, so that $0<r(\theta)<M_0$ for some $M_0>0$.
\medskip

\noindent
{\bf Example} We computed  the curve $\mu=\mu (\xi)$ of $2\pi$-periodic solutions for the problem (\ref{12c}) corresponding to (\ref{be3}), for the case of $a=0.08$ and $b=0.6$, satisfying the condition (\ref{15}). The curve was similar to Figure \ref{fig:1}, with $\mu(0.64) \approx 0$. The polar curve $r=r(\theta)$ centered at $(x_0,y_0)$ (computed at $\xi=0.64$) is then the limit cycle for  Sel'kov's system (\ref{be3}), and it is given in Figure \ref{fig:3}. Again, we used a direct integration to confirm the validity of this limit cycle. 

\begin{figure}
\begin{center}
\scalebox{0.7}{\includegraphics{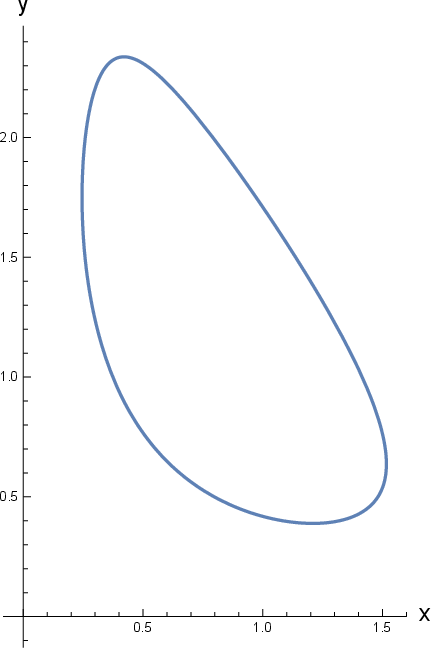}}
\end{center}
\caption{ The limit cycle for Sel'kov's  system (\ref{be3})}
\lbl{fig:3}
\end{figure}
\medskip

However, our computations for Sel'kov's system reveal the limitations of our approach, involving  polar curves around $(x_0,y_0)$: for some $(a,b)$, satisfying the condition (\ref{15}), the solution curve  $\mu=\mu (\xi)$  blows up, before reaching its first root (e.g., for $a=0.001$, $b=0.3$). 
\medskip

\begin{figure}[ht!]

\subfigure[]{\includegraphics[scale=0.8]{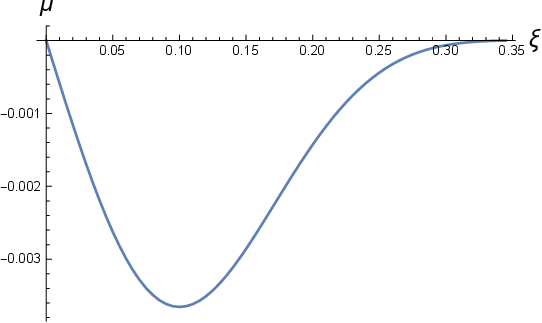}}
\hfill
\subfigure[]{\includegraphics[scale=0.54]{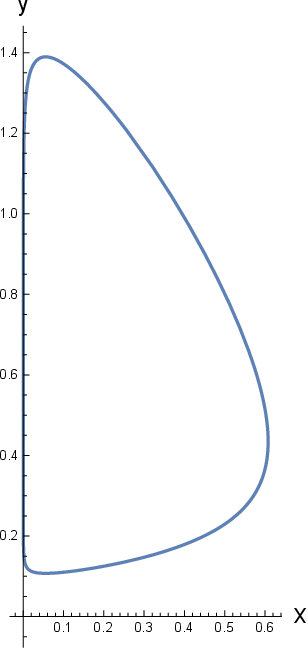}}
\caption{ a)\ The curve of $2\pi$-periodic solutions for the problem (\ref{12c}) corresponding to (\ref{29a}) b)\ The limit cycle for the predator-prey system (\ref{29a}) } 
\label{fig:6}
\end{figure}

A challenging example is presented by a predator-prey system with  the Holling type II functional response, see S.-B. Hsu and J. Shi \cite{shi},
\beq
\lbl{29a}
\frac{dx}{dt}=x(1-x)-\frac{mxy}{a+x} \,, \s \frac{dy}{dt}=-dv+\frac{mxy}{a+x} \,.
\eeq
We computed  the curve $\mu=\mu (\xi)$ of $2\pi$-periodic solutions for the problem (\ref{12c}) corresponding to (\ref{29a}), for the  parameter values $a = 0.5$, 
$m = 1$, and 
$d = 0.1$. This time $\mu(\xi)<0$ for small $\xi>0$, and the first root of $\mu(\xi)$ occurs at $\xi _0 \approx 0.345$. Our computations indicate that $\mu (\xi _0)=\mu' (\xi _0)=0$, and that the curve $\mu=\mu (\xi)$ blows up soon after $\xi _0$. The solution curve for the corresponding equation (\ref{12c}), as well as the   limit cycle of (\ref{29a}), involving a sharp transition,  are given in Figure \ref{fig:6}.

\medskip

A similar approach allows one to  compute multiple limit cycles. For a modified Van der Pol equation
\beq
\lbl{30}
x''+0.04 (x^2-1)(9-x^2)x'+x =0 
\eeq
the energy is being pumped up for small $x$, dissipated at larger $x$, and pumped up again for large $x$. It is natural to expect the existence of a smaller stable limit cycle, and of larger unstable limit cycle. However, the solution curve $\mu=\mu (\xi)$ of the corresponding equation (\ref{12c}) blows up soon after its first root, with $|\mu| \ra \infty$. We see from (\ref{11a}) and  (\ref{17-a}) that  the blow up happens when $\frac{d \theta}{dt}$ approaches zero. Near the second (and any other) limit cycles we expect $\frac{d \theta}{dt}$ to stay away from zero, since otherwise the period of the limit cycle, $T=\int_0^{2\pi} \frac{dt}{d \theta} \, d \theta$, would get large. Hence we need to ``help" the solution curve $\mu=\mu (\xi)$  to get through a ``rough spot", on its way toward its second root. Therefore we suggest approximating the  corresponding equation   (\ref{12c}) by
\beq
\lbl{33}
\frac{dr}{d \theta}+\frac{g(\theta,r)}{1+\epsilon g^{2m}(\theta,r)}=\bar \mu \,,
\eeq	
with some integer $m \geq 1$, and small $\epsilon>0$. For sufficiently small $\epsilon$ the periodic solutions of (\ref{33}) are close to those of  (\ref{12c}). For the example (\ref{30}) above we took $m=3$ and $\epsilon=0.0001$, and the solution curve of (\ref{33}) is presented in Figure \ref{fig:4}. As above, the roots
of $\bar \mu (\xi)$ provided us with very accurate approximation of the limit cycles of the original equation (\ref{30}), see Figure \ref{fig:5}. We verified their accuracy by a direct integration (forward in time  for the smaller stable cycle, and backward in time for the unstable larger one). 

\begin{figure}
\begin{center}
\scalebox{0.7}{\includegraphics{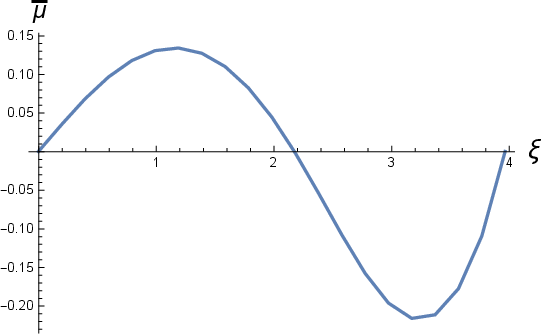}}
\end{center}
\caption{ The curve of $2\pi$-periodic solutions of  (\ref{33})}
\lbl{fig:4}
\end{figure}
\medskip

\begin{figure}
\begin{center}
\scalebox{0.7}{\includegraphics{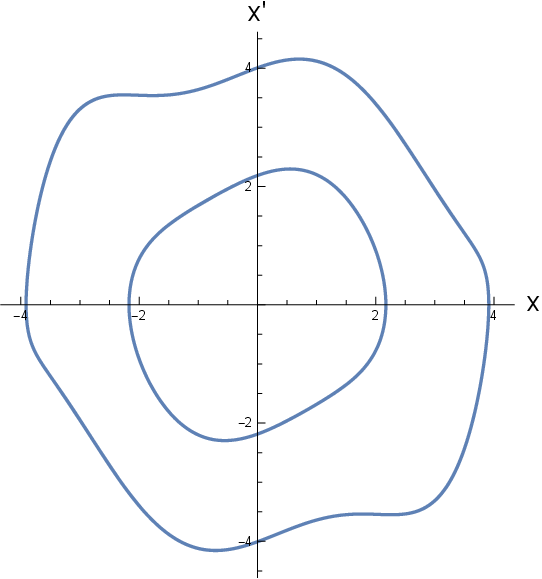}}
\end{center}
\caption{ The two limit cycles   of (\ref{30})}
\lbl{fig:5}
\end{figure}
\medskip

\noindent
{\bf Conclusions}. Given that the autonomous system (\ref{10}) has a rest point at $(x_0,y_0)$, we introduce polar coordinates  with respect to $(x_0,y_0)$, and convert  (\ref{10}) to the periodic first order equation (\ref{12}). Then we embed  (\ref{12}) into the family of equations (\ref{12c}).
Since the functions $g(\theta,r)$ is often singular, we approximate the solutions of  (\ref{12c}) by those of (\ref{33}), for which we compute the global solution curve $\bar \mu=\bar \mu (\xi)$, where $\xi$ is the average of the solution $r(\theta)$ over $[0,2\pi]$. Each root of $\bar \mu (\xi)$ provides us with a limit cycle of (\ref{10}). This straightforward approach avoids the need to construct the trapping region for each particular  system (\ref{10}).
Our approach raises an interesting open question. Suppose a limit cycle of (\ref{10}) encloses just one rest point. Under what condition this  limit cycle
is star-shaped with respect to the rest point (so that it is a polar curve  with respect to the rest point).

\section{Numerical computation of solutions}
\setcounter{equation}{0}
\setcounter{thm}{0}
\setcounter{lma}{0}

To  solve the periodic problem (with $g(t+T,u)=g(t,u)$, for all $t$ and $u$)
\beq
\lbl{n1}
u'+g(t,u)=\mu+e(t) \,, \; u(t)=u(t+T) \;\; \mbox{ for all $t$}\,,
\eeq
we used continuation in $\xi$. (Recall that $u(t)=\xi +U(t)$, with $\int _0^T U(t) \, dt=0$.)
\medskip

\noindent
{\em Step 1}.  We began by implementing the numerical solution of the following linear periodic problem: given the $T$-periodic functions  $b(t)$ and $f(t)$,  find the $T$-periodic solution of
\beq
\lbl{n2}
\s\s\s L[y] \equiv y'(t) +b(t)y=f(t), \;\; y(t)=y(t+T) \,.
\eeq
The general solution of (\ref{n2}) can be written as
\[
y(t)= y_1(t)+c y_2(t) \,,
\]
where $y_1(t)$ is the solution of 
\[
 y'(t) +b(t)y=f(t)\,, \s y(0)=0 \,,
\]
$y_2(t)$ is the solution of 
\[
 y'(t) +b(t)y=0\,, \s y(0)=1 \,,
\]
and $c$ is an arbitrary constant. We have $y(0)=c$,  $y(T)=y_1(T)+c y_2(T)$, and setting $y(0)=y(T)$ gives $c=-\frac{y_1(T)}{y_2(T)-1}$, so that 
\beq
\lbl{n5}
y(t)=y_1(t)-\frac{y_1(T)}{y_2(T)-1} y_2(t)
\eeq
 gives us the $T$-periodic solution of (\ref{n2}), or $L^{-1}[f(t)]$, where $L[y]$ denotes the left hand side of (\ref{n2}), subject to the periodic boundary condition. We used {\em Mathematica}'s NDSolve command to compute $y_1(t)$ and  $y_2(t)$.
\medskip

A natural question is why we did not use the explicit formula available for the equation (\ref{n2}). To solve the nonlinear problem (\ref{n1}) we shall be repeatedly solving the linear problem (\ref{n2}) with $b(t)$ and $f(t)$ given as interpolating functions (produced by {\em Mathematica} on the previous steps). The integrating factor $\nu (t)=e^{\int_0^t b(t) \,dt}$ will require computing a definite integral for each value of $t$, and then solving of (\ref{n2}) will  require double integration, while computation of $y_1(t)$ and $y_2(t)$ can be compared to a single integration. The explicit formula appears to be a useless curiosity, as far as numerical computations go. 
\medskip

\noindent
{\em Step 2}. Then we  implemented the {\em ``linear solver"}, i.e., the numerical solution of the following problem: given the $T$-periodic functions  $b(t)$ and $f(t)$, find the constant $\mu$, so that the problem
\beq
\lbl{n10a}
 y' +b(t)y=\mu+f(t), \; \int_0^T y(t) \, dt=0 
\eeq
has a $T$-periodic solution (of zero average), and compute that solution $y(t)$. The solution is 
\[
y(t)=\mu L^{-1}[1]+L^{-1}[f(t)] \,,
\]
with the constant  $\mu$ chosen so that $\int_0^T y(t) \, dt=0$.
\medskip

\noindent
{\em Step 3}. Turning to the problem (\ref{n1}), we begin with an initial $\xi _0$, and using a step size $\Delta \xi$, we compute the solution $(\mu _n ,u_n(t))$ with the average of  $u_n(t) $ equal to  $\xi _n=\xi _0 +n \Delta \xi$, $n=1,2, \ldots, nsteps$, in the form $u_n(t)=\xi _n+U_n(t)$, where $U_n(t)$ is the $T$-periodic solution of
\beq
\lbl{n3}
 U'+g(t,\xi _n+U)=\mu+e(t) \,, \; \int_0^T U(t) \, dt=0 \,,
\eeq
and $\mu _n$ is the corresponding value of $\mu$. We use Newton's method to approximate the solution of (\ref{n3}) (i.e., $ \left( \mu _n, U_{n}(t) \right)$) by constructing a sequence of approximations $\{ V_m \}$ (with the corresponding values of $\mu$). We set $V_0=U_{n-1}(t)$, and with $V_m$ already computed, 
we linearize the equation (\ref{n3}) at $\xi_{n}+V_{m}(t)$, writing $g(t,\xi _n+U)\approx g(t,\xi _n+V_{m})+g_u(t,\xi _{n}+V_{m})(U-V_{m})$, and call  the linear solver to find the $T$-periodic solution of the problem (\ref{n10a}) with $b(t)=g_u(t,\xi _n+V_{m})$, and $f(t)=e(t)-g(t,\xi _{n}+V_{m})+g_u(t,\xi _{n}+V_{m})V_{m}$, obtaining an approximation of $V_m$ and of the corresponding $\mu$.
We found that just four iterations of Newton's method ($1 \leq m \leq 4$), coupled with a relatively small $\Delta \xi$ (e.g., $\Delta \xi=0.1$), were sufficient for accurate computation of the solution curves. Finally, we plot the points $(\xi _n, \mu _n)$ to obtain the solution curve.
\medskip

We  verified our numerical results by an independent calculation. Once  a periodic solution $u(t)$ of (\ref{n1})  is computed at some $\mu$, we took its value of $u(0)$, and computed numerically the solution of (\ref{n1}) at this $\mu$, with the initial value equal to $u(0)$ (using the NDSolve command), as well as the  average of $u(t)$. We had a complete agreement (with $u(t)$ and $\xi$) for all $\mu$, and all equations that we tried.
\medskip

For applications to fishing we need to solve the periodic problem (with $g(t+T,u)=g(t,u)$, for all $t$ and $u$)
\beq
\lbl{n20}
u'+g(t,u)=\mu j(t) \,, \; u(t)=u(t+T) \;\; \mbox{ for all $t$} \,,
\eeq
where $j(t)$ is  $T$-periodic, and $j(t)<0$.
Similarly to the above, we used continuation in $\xi$. The Step 1 is as above. The ``linear solver" part is modified:  given the $T$-periodic functions  $b(t)$, $j(t)$ and $f(t)$, find the constant $\mu$, so that the problem
\beq
\lbl{n10}
 y' +b(t)y=\mu j(t)+f(t), \; \int_0^T y(t) \, dt=0 
\eeq
has a $T$-periodic solution, and compute that solution $y(t)$. The solution is 
\[
y(t)=\mu L^{-1}[j(t)]+L^{-1}[f(t)] \,,
\]
with the constant  $\mu$ chosen so that $\int_0^T y(t) \, dt=0$. 
\medskip

To solve (\ref{n20}), we begin with an initial $\xi _0$, and using a step size $\Delta \xi$, we compute the solution $(\mu _n ,u_n(t))$ with the average of  $u_n(t) $ equal to  $\xi _n=\xi _0 +n \Delta \xi$, $n=1,2, \ldots, nsteps$, in the form $u_n(t)=\xi _n+U_n(t)$, where $U_n(t)$ is the $T$-periodic solution of
\beq
\lbl{n23}
 U'+g(t,\xi _n+U)=\mu j(t) \,, \; \int_0^T U(t) \, dt=0 \,,
\eeq
and $\mu _n$ is the corresponding $\mu$. We use Newton's method to approximate the solution of (\ref{n23}), by constructing a sequence  of approximations $\{ V_m \}$. We set $V_0=U_{n-1}(t)$, and with $V_m$ already computed, 
we linearize the equation (\ref{n3}) at $\xi_{n}+V_{m}(t)$, writing $g(t,\xi _n+U)\approx g(t,\xi _n+V_{m})+g_u(t,\xi _{n}+V_{m})(U-V_{m})$, and call  the new linear solver to find the $T$-periodic solution of the problem (\ref{n10}) with $b(t)=g_u(t,\xi _n+V_{m})$, and $f(t)=-g(t,\xi _{n}+V_{m})+g_u(t,\xi _{n}+V_{m})V_{m}$, obtaining an approximation of $V_m$ and the corresponding $\mu$. After around four steps of Newton's method, we obtain an accurate approximation of $\left( \mu _n, U_n(t) \right)$.
Finally, we plot the points $(\xi _n, \mu _n)$ to obtain the solution curve.

\section{Periodic population model with fishing}
\setcounter{equation}{0}
\setcounter{thm}{0}
\setcounter{lma}{0}

We begin with the logistic population model ($u=u(t)$ is the number of units of fish at time $t$)
\[
u'=u(a(t)-u) \,.
\]
We assume that $a(t)$ is a continuous and $T$-periodic function, satisfying $\int_0^T a(t)\,dt>0$. By direct integration, it is easy to see that there is a unique $T$-periodic solution, which attracts all other positive solutions, as $t \ra \infty$, see e.g.,  M.N. Nkashama \cite{N}, or  P. Korman \cite{K}  (exact multiplicity results for more general polynomial models were given in P. Korman and T. Ouyang \cite{KO}). To model fishing, let $j(t)$ be a  positive continuous $T$-periodic function, giving  the number of units of fish taken at time $t$, and let the number $b>0$ give the decrease in the growth rate per one unit of fish taken, and denote $f(t)=bj(t)$. Then the model is
\[
u'=u(a(t)-u)-f(t) \,, \s\s u(t+T)=u(t) \s \mbox{for all $t$}\,.
\]
Let us now vary the intensity of fishing, and consider 
\beq
\lbl{p1}
u'=u(a(t)-u)-\mu f(t) \,, \s\s u(t+T)=u(t) \s \mbox{for all $t$}\,,
\eeq
depending on a parameter $\mu>0$.
\medskip

The following result  is well known see \cite{A}, \cite{M}, \cite{N}, and also \cite{AL}. We sketch the proof for completeness.
\begin{prop}\lbl{prop:1}
Assume that  $\int_0^T a(t)\,dt>0$, and $f(t)>0$.  Then there is a critical $\mu _0 >0$, so that the problem (\ref{p1}) has two exactly positive $T$-periodic solutions for $0<\mu<\mu _0$, exactly one positive $T$-periodic solution at $\mu=\mu _0$ and for $\mu <0$, and no $T$-periodic solution for $\mu>\mu _0$. (With the bifurcation diagram as in Figure \ref{fig:2}. The case when  $\mu <0$ can be interpreted as  stocking of fish). Moreover, in case $\mu>\mu _0$, all solutions of (\ref{p1}) tend to $-\infty$, in finite time.
\end{prop}

\pf
At $\mu =0$ the equation (\ref{p1}) has a unique $T$-periodic solution. We now continue $T$-periodic solutions for $\mu >0$. Clearly, $T$-periodic solutions remain positive. Indeed, assuming otherwise, we can find an interval $(t_1,t_2)$ on which $u(t)<0$ that $u( t_1)=u( t_2)=0$. From the equation (\ref{p1}), $u'(t)<0$ on $(t_1,t_2)$, a contradiction. The continuation in $\mu$ cannot proceed indefinitely, because the right hand side of  (\ref{p1}) becomes negative. There is at most two positive $T$-periodic solutions of  (\ref{p1}), see e.g. p. 245 in \cite{K}.
\epf  

It follows from the Proposition  \ref{prop:1} that  for $\mu>\mu _0$ the fishing is not sustainable, while $\mu=\mu _0$ gives the maximal  sustainable level of fishing. We call  $\mu _0 \int_0^T f(t) \, dt$ {\em the maximal fishing take} ($\frac{\mu _0}{T}  \int_0^T j(t) \, dt$ gives the maximal possible average number of fish, caught over one  period). The problem is:
find the fishing strategy $f(t)$ (or  $j(t)$) to maximize the maximal fishing take  $ \mu _0 \int_0^T f(t) \, dt $. For example, assume that the growth rate $a(t)$ has one point of  minimum and one point of   maximum, corresponding to summer and winter. What are good seasons for fishing?
\medskip

In case $a(t)=a$ and $f(t)=1$ are constants, the  $T$-periodic solutions of (\ref{p1}) are also constants. In that case $\mu _0=\frac{a^2}{4}$, and the maximal fishing take is then $ \mu _0 T=\frac{a^2}{4} T $. 

\begin{prop}
Assume that $a(t)=a$ is constant. Then the optimal fishing strategy is for $f(t)$ to be a constant too, i.e.,
\beq
\lbl{p5}
\mu _0 \int_0^T f(t) \, dt \leq \frac{a^2}{4} T \,.
\eeq
\end{prop}

\pf
Assume that $u_0(t)$ is the solution of (\ref{p1}) at $\mu _0$, i.e.,
\beq
\lbl{p2}
u'_0=au_0-u_0^2-\mu _0 f(t)\,.
\eeq
Since the solution curve does not continue in $\mu$ to the right of $\mu _0$, the implicit function theorem does not apply, and so the corresponding linearized equation 
\[
w'=(a-2u_0)w
\]
has non-trivial $T$-periodic solutions. By Lemma \ref{lma:g2}, $\int _0^T u_0(t) \, dt=\frac12 aT$. We now integrate (\ref{p2}), to obtain
\[
\mu _0 \int_0^T f(t) \, dt=\frac12 a^2T-\int _0^T u^2_0(t) \, dt \leq \frac12 a^2T-\frac{1}{T} \left( \int _0^T u_0(t) \, dt \right)^2=\frac{a^2}{4} T \,,
\]
concluding the proof.
\epf

The inequality (\ref{p5}) can be generalized to any $T$-periodic $a(t)$. Indeed writing the equation (\ref{p1}) as 
\[
u'=-\left(u-\frac12 a(t) \right)^2+\frac14 a^2(t)-\mu f(t) \,,
\]
and integrating, we conclude that for all $\mu$  
\beq
\lbl{p6}
\mu  \int_0^T f(t) \, dt \leq \frac14 \int_0^T  a^2(t) \, dt \,.
\eeq

Suppose  that $a(t)$ is not a constant. Is it true that the optimal fishing strategy is still for $f(t)$ to be a constant? Our  numerical experiments indicate that the answer is no, in general.

\begin{figure}
\begin{center}
\scalebox{0.75}{\includegraphics{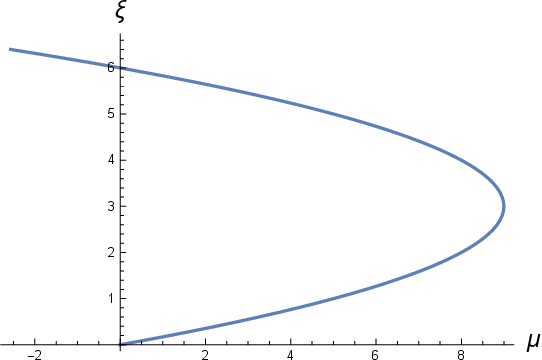}}
\end{center}
\caption{ The curve of $1$-periodic solutions of (\ref{p10})}
\lbl{fig:2}
\end{figure}

\medskip

\noindent
{\bf Example} We computed the curve of $1$-periodic solutions for
\beq
\lbl{p10}
u'=\left(1+0.4 \sin 2 \pi t \right)u-\mu \left(1+b \sin 2 \pi t \right) \,,
\eeq
with $b=0.2$, see Figure \ref{fig:2}. Here the maximal fishing take turned out to be  $\mu _0 \int_0^1 \left(1+b \sin 2 \pi t \right) \, dt  =\mu _0 \approx 8.9955$. When we took $b=0$, the maximal fishing take decreased to $ \approx 8.9818$. In Figure \ref{fig:2}, $\xi=\int_0^1 u(t) \, dt$, and our computations were using the method developed in the preceding section.
\medskip

\noindent
{\em Open Problem:} Given $T$-periodic non-constant $a(t)$, find the optimal  $T$-periodic $f(t)$, maximizing the maximal fishing take $\mu _0 \int_0^T f(t) \, dt$.
\medskip

Our numerical experiments  suggest the following conjecture. 
\medskip

\noindent
{\em Conjecture}: Decompose $a(t)=A+\al(t)$, with constant $A>0$ and $\int _0^T \al (t) \, dt=0$. Then  the maximal  fishing take satisfies
\beq
\lbl{p7}
\mu _0 \int_0^T f(t) \, dt \leq \frac{A^2}{4} T \,.
\eeq
Observe that $\frac14 \int_0^T  a^2(t) \, dt=\frac{A^2}{4} T+\frac14 \int_0^T  \al ^2(t) \, dt$, so that the inequality (\ref{p7}) is significantly sharper than  (\ref{p6}). 
\medskip

The inequality (\ref{p7}) implies that if the averages of $a(t)$ and $f(t)$ are fixed, the best fishing outcome occurs if $a(t)$ and $f(t)$ are constants.

\end{document}